\documentclass[11pt,a4paper,reqno]{amsart}
\usepackage{multirow}
\usepackage[centertags]{amsmath}
\usepackage{amsfonts}
\usepackage{amssymb}
\usepackage{amsthm}
\usepackage{newlfont}
\usepackage{mathtools}
\usepackage{a4}
\usepackage{enumerate}
\usepackage[pdftex]{graphicx,color}
\usepackage{diagbox}

\newcommand{\ben}{\begin{equation}}
\newcommand{\een}{\end{equation}}
\newcommand{\bens}{\begin{equation*}}
\newcommand{\eens}{\end{equation*}}
\newcommand{\bal}{\begin{align}}
\newcommand{\eal}{\end{align}}
\newcommand{\nn}{\nonumber\\ }

\newcommand{\cB}{{\cal B}}

\newcommand{\ch}{{\text{ch}}}

\newcommand{\mult}{{\text{mult}}}
\newcommand{\hfb}{{\hfill\break}}

\input amssym.def
\def\Z{{\mathbb Z}} \def\R{{\mathbb R}}  \def\N{{\mathbb N}}
   
\begin{document}
\baselineskip=16pt
\parskip=10pt

\begin{center}
{\LARGE {\bf Demazure formula for \scalebox{1.15}{$A_{\lowercase{n}}$} Weyl polytope sums}}

\vspace{6mm}
{\LARGE J{\o}rgen Rasmussen$^1$ and Mark A. Walton$^2$}
\\[.4cm]
{\em $^1$School of Mathematics and Physics, University of Queensland}\\
{\em St Lucia, Brisbane, Queensland 4072, Australia}
\\[.2cm]
{\em $^2$Department of Physics \& Astronomy, University of Lethbridge}\\
{\em Lethbridge, Alberta\  T1K 3M4, Canada}
\\[.4cm]
{\tt j.rasmussen\,@\,uq.edu.au}\qquad {\tt walton\,@\,uleth.ca}
\end{center}

\vspace{12mm}
\centerline{{\bf{Abstract}}}
\vskip.4cm
\noindent
The weights of finite-dimensional representations of simple Lie algebras are naturally associated with Weyl polytopes. Representation characters decompose into multiplicity-free sums over the weights in Weyl polytopes. The Brion formula for these Weyl polytope sums is remarkably similar to the Weyl character formula. Moreover, the same Lie characters are also expressible as Demazure character formulas. This motivates a search for new expressions for Weyl polytope sums, and we prove such a formula involving Demazure operators. It applies to the Weyl polytope sums of the simple Lie algebras $A_n$, for all dominant integrable highest weights and all ranks $n$.

\vfill\eject

\section{Introduction} 

The Brion formula \cite{{Brion88}, {Brion92}} is a general expression for an  exponential sum over lattice points in 
a polytope.  That sum is sometimes called the integer-point transform of the polytope. 

In the weight lattice of a simple Lie algebra, a Weyl polytope has the weights in an orbit of the 
Weyl group as its vertices.  Applied to Weyl polytopes, the Brion theorem yields a formula that is remarkably similar to the 
Weyl character formula \cite{DhillonKhare16, Walton04, Postnikov09}.  
As a consequence, the {\it polytope expansion} of Weyl characters in terms of integer-point transforms is 
natural and useful \cite{Walton04, Postnikov09, Schutzer12, Walton13, DhillonKhare16, Rasmussen18}.  

Here we explore further the relation between Lie characters and integer-point transforms.  
Other formulas exist for the characters, and following \cite{Walton04}, 
we are interested in finding similar expressions for the  
integer-point transforms of Weyl polytopes, herein called {\it Weyl polytope sums}.
We focus on the Demazure character formula \cite{Demazure74, Andersen85, Joseph85} and use the 
Demazure operators involved to write expressions for the Weyl polytope sums.   
We have obtained results for the simple Lie algebras $A_n$, for all ranks $n\in\mathbb{N}$.  

In the following section, we review the initial motivation for the present work, in part to establish our notation. 
We describe the similarity between the Weyl 
character formula and  the Brion formula, and the polytope expansion that exploits it. Section 3 
is a quick account of the Demazure character formula.  Our new formula for $A_n$ Weyl polytope sums is 
presented and proved in Section 4.  The final section offers a short conclusion. 

\section{Polytope expansion of Lie characters}
  
Let $X_n$ denote a simple Lie algebra of rank $n$ ($X$ is a letter from $A$ to $G$). The sets of 
fundamental weights and simple roots 
are denoted by $F:=\{\Lambda^i\,\vert\, i=1,\ldots,n\}$ and $S:=\{\alpha_i\,\vert\, i=1,\ldots,n\}$, 
respectively. The corresponding weight and root lattices are $P:=\Z\,F$ and $Q:= \Z\, S$.
The set of dominant integrable weights is $P_+ := \N_0 F$, and we write $R$ 
($R_+$, $R_-$) for the set of (positive, negative) roots of $X_n$. 

\subsection{Weyl character formula} 

Consider a finite-dimensional irreducible module $L(\lambda)$ over $X_n$ of highest weight $\lambda\in P_+$. 
The formal character of $L(\lambda)$ is defined as
\begin{align}
\ch_\lambda\ :=\ \sum_{\mu\in P}\, \mult_\lambda(\mu)\, e^\mu
\ =\ \sum_{\mu\in P(\lambda)}\, \mult_\lambda(\mu)\, e^\mu\, ,
\label{char}\end{align}
where $\mult_\lambda(\mu)$ is the multiplicity of the weight $\mu$ in the module $L(\lambda)$, 
while $P(\lambda)$ is the set of weights of $L(\lambda)$:
\begin{align}
P(\lambda)\ =\ \{\, \mu\in P \,\vert\, \mult_\lambda(\mu)\, >\, 0 \,\}\ \ .
\end{align}

The formal exponentials of weights obey $e^\mu\, e^\nu=e^{\mu + \nu}$. With
\begin{align}
e^\mu(\sigma)\ :=\ e^{( \mu, \sigma )}\  , \qquad \sigma\in P\ ,
\label{FormalExp}
\end{align}
where $( \mu, \sigma )$ is the inner product of  weights $\mu$ and $\sigma$, the formal exponential 
$e^\mu$ simply stands for $e^{( \mu, \sigma )}$ before a choice of weight $\sigma$ is made.  A choice 
of $\sigma$ fixes a conjugacy class of elements in the Lie group $\exp(X_n)$, and the formal character becomes 
a true character: $\ch_\lambda(\sigma)$, the trace, in the irreducible highest-weight
representation of highest weight $\lambda$, of elements 
of $\exp(X_n)$ in the conjugacy class labelled by  $\sigma$. 
 
The celebrated Weyl character formula is 
\begin{align}
\ch_\lambda\ =\ \frac{\sum_{w\in W} (\det w)\, e^{w(\lambda+\rho)-\rho}}{\prod_{\alpha\in R_+}\, (1-e^{-\alpha})}\ , 
\label{WCF}\end{align} 
where $\rho:=\sum_{i=1}^n{\Lambda^i}$.
The Weyl invariance of the character can be made manifest, as
 \begin{equation}
\ch_\lambda\ =\
\sum_{w\in W}\, e^{w\lambda}\, \prod_{\alpha\in  R_+}\,
(1-e^{-w\alpha})^{-1}\ ,
\label{WCFformal}\end{equation}
where
\begin{equation}
(1-e^\beta)^{-1}\ = \left\{\begin{matrix} 1+e^\beta+e^{2\beta}+\ldots\,, &\quad\ \beta\in R_-\ ; 
\\[.2cm]  
-(e^{-\beta}+e^{-2\beta}+\ldots)\,, &\quad\ \beta\in R_+\ . \end{matrix} \right. 
\label{series}\end{equation}
For brevity, we sometimes write $G_\beta:=(1-e^{-\beta})^{-1}$, for $\beta\in R$.
For each $w\in W$, we also define 
\begin{align} w(e^\mu)\, :=\, e^{w\mu}\ ,\qquad w\, e^\mu\,:=\, e^{w\mu}\, w\ ,
\label{wop}
\end{align}
meaning that $w$ acts on an explicitly indicated argument only or, if no argument is given, on everything to its right.
We can then rewrite (\ref{WCFformal}) as 
\begin{equation}
\ch_\lambda\ =\ {\cal C}(e^\lambda)\ ,
\label{WCFop}
\end{equation}
with
\begin{equation}
{\cal C}\ :=\
\sum_{w\in W} \Big(\prod_{\alpha\in  R_+}\,
G_{w\alpha}\Big)\,w\ =\ \sum_{w\in W}\,w\,\prod_{\alpha\in  R_+}\,
G_{\alpha}\ . 
\label{WCFopC}\end{equation}

\subsection{Brion formula}

A polytope is the convex hull of finitely many points in $\R^d$. A polytope's vertices form such a set of points, 
with minimum cardinality. A lattice polytope has all its vertices in an integral lattice in $\R^d$.  
The corresponding (formal) integer-point transform of the polytope is the sum of terms $e^\phi$ 
over the lattice points $\phi$ in the polytope. 

Brion \cite{Brion88, Brion92} found a general formula for these integer-point transforms. For $\lambda\in P$, 
let the {\textit{Weyl polytope}} ${\text Pt}_\lambda$ be the polytope with vertices given by the Weyl orbit 
$W\lambda$.  Consider the integer-point transform  
\begin{equation}   
{\text B}_\lambda\ :=\ \sum_{\mu\in (\lambda+Q)\cap{\text Pt}_\lambda}\, e^\mu\ ,
\label{Bdefinition}
\end{equation}
where the relevant lattice is the $\lambda$-shifted root lattice $\lambda + Q$ of the algebra $X_n$.  
We refer to these integer-point transforms as {\it Weyl polytope sums}.  Since 
\begin{align}
{\text B}_\lambda\ =\ \sum_{\mu\in P(\lambda)}\,  e^\mu\, ,
\label{Bmfchar}\end{align} 
the Weyl polytope sum has an interpretation as a ``multiplicity-free'' character \cite{DhillonKhare16}.
By (\ref{Bmfchar}), it is obtained from the character (\ref{char}) by putting $\mult_\lambda(\mu) \to 1$ 
for all $\mu\in P(\lambda)$. 

Applied to a 
Weyl polytope, the Brion formula yields
\begin{equation}
{\text B}_\lambda\ =\
\sum_{w\in W}\, e^{w\lambda}\, \prod_{{\alpha\in }S}\,
(1-e^{-w\alpha})^{-1}\ . 
\label{Brion}\end{equation}
Following (\ref{WCFop}) and (\ref{WCFopC}), we rewrite (\ref{Brion}) as 
\begin{equation}
{\text B}_\lambda\ =\ {\cal B}(e^\lambda)\ ,
\label{Brionop}
\end{equation}
with
\begin{equation}
{\cal B}\ :=\
\sum_{w\in W}\Big( \prod_{\alpha\in  S}\,
G_{w\alpha}\Big)\,w\ =\ \sum_{w\in W}\, w\, \prod_{\alpha\in  S}\,
G_{\alpha}\ . 
\label{BrionopB}\end{equation}

\subsection{Polytope expansion}
 
The Brion formula (\ref{Brion}) is remarkably similar to the Weyl character formula
(\ref{WCFformal}) \cite{Walton04, Postnikov09, DhillonKhare16}.  It is therefore natural and fruitful to consider 
the polytope expansion of Lie characters \cite{Walton04, DhillonKhare16, Walton13}:
\begin{equation}  {\text{ch}}_\lambda\ =\ \sum_{\mu\leq\lambda}\, A_{\lambda,\mu}\,\, {\text B}_\mu\ .\ 
\label{PolytExpn}\end{equation} 
The expansion coefficients $A_{\lambda,\mu}$ are integers.  They were dubbed  {\em polytope multiplicities} and denoted ${\text{polyt}}_\lambda(\mu)$ in \cite{Walton13}, in analogy 
with the weight multiplicities ${\text{mult}}_\lambda(\mu)$ appearing in the expansion (\ref{char}).   For type $A$, they were shown in \cite{LecouveyLenart21}
to be non-negative.  However, other examples have been found that are negative  \cite{LecouveyLenart21}, so ``multiplicity" appears to be a misnomer.

We do not consider the polytope expansion further in this note. 
Instead, we focus on the striking relationship between characters and Weyl polytope sums.

\section{Demazure character formula}

Here we show that expressions similar to the Demazure character formula can be written for the 
Weyl polytope sums in (\ref{Bdefinition}).

Let us first sketch the Demazure character formula. The Weyl group $W$ is generated by the reflections 
$r_\beta$  in weight space across the hyperplanes normal to the corresponding roots $\beta\in R$:  
\begin{align}
r_\beta(\lambda)\ :=\ \lambda\, -\, (\lambda,\beta^\vee)\,\beta\ ,
\label{PrRefl}
\end{align}
where $\beta^\vee:= 2\beta/(\beta,\beta)$. 
In fact, the Weyl group is generated by the primitive (simple-root) reflections, $r_i\equiv r_{\alpha_i}$.

For each primitive reflection $r_i$, we define the Demazure operator
\begin{align}
D_i\ :=\ \frac{1\ -\ e^{-\alpha_i}\,  r_i}{1-e^{-\alpha_i}}\ .
\label{Diri}
\end{align}
For $\lambda\in P$, we set $r_i(e^\lambda):=e^{r_i\lambda}$ and thus have
\begin{align}
D_i(e^\lambda)\ =\ \left\{\begin{array}{ll}  
 e^{\lambda}+ e^{\lambda-\alpha_i} + e^{\lambda-2\alpha_i} + \ldots + e^{r_i\lambda}\, ,
 \quad\ & (\lambda,\alpha_i^\vee)\ge 0\ ; 
 \\[0.2cm] 
 0\, , \quad\ &(\lambda,\alpha_i^\vee) = -1\ ; 
\\[0.2cm] 
 -\,(e^{\lambda+\alpha_i} + e^{\lambda+2\alpha_i} + \ldots + e^{r_i(\lambda+\alpha_i)})\,, 
 \quad\ & (\lambda,\alpha_i^\vee)  < -1\ .\end{array}\right.
\label{Di}
\end{align}
We will also use the {\it modified Demazure operators} 
\begin{align}
d_i\ :=\ D_i\,-\,1\ =\ \frac{e^{-\alpha_i}\,\left( 1\, -\,   r_i\right)}{1-e^{-\alpha_i}}\ . 
\label{didef}
\end{align}

For every $w\in W$, a Demazure operator $D_w$ can be defined:
In a reduced decomposition of $w$, replace the factors $r_j$ with $D_j$.  
Demazure has shown \cite{Demazure74} that the resulting operator $D_w$ is  independent of which reduced decomposition is 
used (see also \cite{Joseph85}). Accordingly, the Demazure operators obey relations encoded in the Coxeter-Dynkin diagrams of $X_n$. 
To illustrate, let $w_L\in W$ denote the longest element of the Weyl group.
For $A_2$, for example, we have
$w_L=r_1r_2r_1=r_2r_1r_2$, and the associated Demazure operator can be written in the two ways
$D_{w_L}=D_1D_2D_1=D_2D_1D_2.$

The Demazure character formula takes the form (\ref{WCF}) with
\begin{equation}  
{\cal C}\ =\  D_{w_L}\ ,\qquad {\text i.e.\ }\qquad \ch_\lambda\ =\ D_{w_L}(e^{\lambda})\ .
\label{DemazureD}\end{equation}

\section{$A_n$ weight-polytope formulas of Demazure type}

We will now restrict attention to the simple Lie algebras $A_n$.
When appropriate, the superscript $(n)$ will be used to indicate the dependence on the rank $n$. 

With $G_i\equiv G_{\alpha_i}$, rewriting
\begin{align}
D_i\ =\ (1+r_i)\,G_i
\label{DwithG}\end{align} 
will be useful. For $k,m\in\{1,\ldots,n\}$ with $k\neq m$, we also define
\begin{align}
G_{k,m}\ :=\ (1-e^{-\alpha_k-\alpha_m})^{-1}\ .
\end{align}
From
\begin{align}
G_k\, r_m\ =\ (1-e^{-\alpha_k})^{-1}\, r_m\ =\ r_m\, (1-e^{-r_m\alpha_k})^{-1}\ ,
\label{GrrG}\end{align}
it follows that 
\begin{align}
G_{k+1}\,r_k\ =\ r_k\, G_{k,k+1}\, ,\ \ &\quad G_{k-1}\,r_k\ =\ r_k\, G_{k-1,k}\ , \nn 
G_{k,k+1}\,r_k\ =\ r_k\, G_{k+1}\, ,\ \ &\quad G_{k-1,k}\,r_k\ =\ r_k\, G_{k-1}\ .
\label{GrrGkk}\end{align}
Here and henceforth we use the standard numbering of $A_n$ simple roots, so that $(\alpha_k, \alpha^\vee_{k+1})=-1$ for all $k\in \{1,\ldots,n-1\}$. 

Let us define 
\begin{align}
s_{i,j}\ :=\ r_j r_{j-1}\cdots r_i\ ,\qquad 1\le i\le j\le n\ ,
\label{sij}\end{align}
and
\begin{align}
w_{i,j}\ :=\ \sum_{k=i}^j\, s_{k, j}+1\ =\ s_{i,j}\, +\,  s_{i+1,j} + \ldots + s_{j,j} + 1\ .
\label{wijdef}\end{align}
The following expression is reminiscent of the Poincar\'e series discussed by Macdonald in \cite{Macdonald72}.

\noindent{\bf Lemma.} 
\ For the simple Lie algebras $A_n$, 
\begin{align}
\sum_{w\in W^{(n)} } w\ =\ w_{1,1}\, w_{1,2}\, \cdots\, w_{1,n}\ .
\label{Wnsum}\end{align}
{\bf Proof.}
Since $W^{(n)}\cong S_{n+1}$, the result follows by showing that the sum in (\ref{Wnsum}) acting on 
$(1,\ldots,n+1)$ produces all possible permutations thereof, 
where $r_j(1,\ldots, j-1,j,j+1,j+2,\ldots,n+1):= (1,\ldots, j-1,j+1,j,j+2,\ldots,n+1)$, for $j=1,\ldots,n$.
Our proof is by induction on $n$.
For $n=1$, the result is trivial; for $n=2$, it is readily verified. 
Applying $w_{1,n}$ to $(1,\ldots,n+1)$ produces $n+1$ terms, one for each possible value of the last entry. 
By the induction hypothesis, the relation (\ref{Wnsum}) holds for $n$ replaced by $n-1$,
so all permutations of $(1,\ldots,n+1)$ result. 
$\blacksquare$

\noindent
Our main result is given in (\ref{cBD}) below and may be viewed as the Demazure analogue of this Lemma.

For $1\le i\le j\le n$, we now define the {\em generalized Demazure operators}
\begin{align}
D_{i,j}\ :=\ r_j r_{j-1}\cdots r_{i+1}\,d_i\,+\,r_j r_{j-1}\cdots r_{i+2}\,d_{i+1}\,+ \ldots +\,d_{j}\,+\,1\ ,
\label{Dijdef}\end{align}
noting that $D_{i,i}=d_i+1=D_i$.

\noindent {\bf Lemma.}
For $1\le i\le j\le n$, 
\begin{align}
D_{i,j}\ =\ s_{i,j} G_i&\,+\, s_{i+1,j} G_iG_{i+1}/G_{i,i+1}\, +\, s_{i+2,j}  G_{i+1}G_{i+2}/G_{i+1,i+2}\nn 
&\, +\ldots+\, s_{j,j} G_{j-1}G_{j}/G_{j-1,j}\, +\, G_j \ .
\label{DijGGG}\end{align}
{\bf Proof.}
Use (\ref{DwithG}) and substitute
\begin{align}
d_i \ =\ (r_i+1)G_i\, -\, 1
\label{dwithG}\end{align}
into (\ref{Dijdef}) to obtain
\begin{align}
D_{i,j}\ =\ s_{i,j} G_i&\,+\,s_{i+1,j} (G_i-1+G_{i+1})\,+\,s_{i+2,j}  (G_{i+1}-1+G_{i+2})\nn
&\,+\ldots+\,s_{j,j} (G_{j-1} -1+G_j)\,+\,G_j \ .
\label{DijG}\end{align}
Now notice that 
\begin{align}
G_k\, -\, 1\, +\, G_m\ =\ G_k\, G_m\, (1-e^{-\alpha_k-\alpha_m})\ =\ G_k\, G_m\, /G_{k,m}\ ,
\label{Gkm}\end{align}
from which the result follows. 
$\blacksquare$

With the similarity between (\ref{Dijdef}) and (\ref{wijdef}), 
the Weyl group algebra relation (\ref{Wnsum}) motivates our main result, 
the following Theorem.  

\noindent{\bf Theorem.}
\ For the simple Lie algebras $A_n$, 
\begin{align}
{\cal B}^{(n)}\ =\ D_{1,1}\, D_{1,2}\, \cdots\, D_{1,n}\ . 
\label{cBD}\end{align} 
{\bf Proof.}
Our proof is by induction on the rank $n$.  
For $n=1$, the characters and Weyl polytope sums coincide, and $w_L=r_1$, so 
\begin{align}
{\cal B}^{(1)}\  =\ {\cal C}^{(1)}\ =\ D_{w_L}=\ D_1\ =\ D_{1,1}\  .
\label{none}\end{align} 
For the induction step, we assume 
\begin{align}
D_{1,1}\, D_{1,2}\, \cdots\, D_{1,n-1}\ =\ w_{1,1}\, w_{1,2}\, \cdots\, w_{1,n-1}(G_1G_2\cdots G_{n-1})\ , 
\label{DwGnm}\end{align}
and prove that 
\begin{align}
D_{1,1}\, \cdots\, D_{1,n-1}\, D_{1,n}\ =\ w_{1,1}\,  \cdots\, w_{1,n-1}\, w_{1,n}(G_1\cdots G_{n-1}G_n)\ . 
\label{DwGn}\end{align}
By (\ref{DwGnm}), the left-hand-side of (\ref{DwGn}) is 
\begin{align}
w_{1,1}\, w_{1,2}\, \cdots\, w_{1,n-1}(G_1G_2\cdots G_{n-1})D_{1,n}\ .
\end{align}
To complete the proof, we establish
\begin{align}
(G_1G_2\cdots G_{n-1})D_{1,n}\ =\ w_{1,n}(G_1G_2\cdots G_{n-1}G_n)\ .
\label{GsDwGs}\end{align}
By (\ref{DijGGG}), we have
\begin{align}
D_{1,n}\ =\ s_{1,n} G_1\, +\, s_{2,n} G_1G_{2}/G_{1,2}\, +\, s_{3,n}  G_{2}G_{3}/G_{2,3}
\,+\ldots +\, s_{n,n} G_{n-1}G_{n}/G_{n-1,n}\, +\, G_n\ .
\label{DonenGGG}\end{align}
Using (\ref{GrrGkk}), we see that 
\begin{align}
(G_1 \cdots G_{n-1})\, s_{j,n}\ =\ s_{j,n}\, (G_1\cdots G_{j-2}\, G_{j-1,j}\, G_{j+1}\cdots G_n)\ .
\end{align}
The relation (\ref{GsDwGs}) now follows, thus completing the proof. 
$\blacksquare$

\section{Conclusion} 

Our main result is the formula (\ref{cBD}) involving (modified) Demazure operators for the 
weight-polytope lattice sums for the Lie algebras $A_n$. 
It is valid for all ranks $n\in \N$ and all dominant integrable highest weights. 

In \cite{MWproc}, formulas are written for the rank-2 weight-polytope lattice sums. It is interesting to note that 
these formulas are easily recast into a form that is very similar to the one we have found for $A_n$. 
Apart from $A_2$, $C_2\cong B_2$ and $G_2$ are the only (up to isomorphism) rank-2 simple Lie algebras. 
In both cases, let $\alpha_1$ denote the short root. For $C_2$, we then find
\begin{align}
\cB\ =\ (1+d_2)\,(1+d_1+r_1d_2+r_1r_2d_1)\ ,
\label{Ctwo}
\end{align}
while for $G_2$, we obtain 
\begin{align}
\cB\ =\ (1+d_2)\,(1+d_1+r_1d_2+r_1r_2d_1+r_1r_2r_1d_2+r_1r_2r_1r_2d_1)\ . 
\label{Gtwo}
\end{align}
We believe this indicates that we are on track toward a general form, valid for all simple Lie algebras. 
Furthermore, we hope that such a formula might lead to one that applies beyond the Lie context, 
as the Brion formula does, to polytopes besides the Weyl polytopes.  

To finish, let us mention some interesting related work. The polytope expansion of Lie characters is highly 
reminiscent of the early work of Antoine and Speiser \cite{AntoineSpeiser64} and the 
recursive formulas found by Kass \cite{Kass91}. Recent work generalizes the context significantly.
Dhillon and Khare \cite{DhillonKhare16} thus report results for all simple highest-weight
modules over Kac-Moody algebras.  In \cite{LecouveyLenart21} by Lecouvey and Lenart, a connection with 
the {\it atomic decomposition} of characters is described, along with ($q$- or $t$-)deformations of the 
structures described herein.

\section*{{Acknowledgements}}
\vskip-0.2cm
The research of J.R. was supported by the Australian Research Council under the Discovery Project scheme, 
project number DP200102316. 
The research of M.W. was supported by a Discovery Grant (40082-DG) from the Natural Sciences and Engineering 
Research Council of Canada (NSERC).
The authors thank Ole Warnaar for comments.

\bibliographystyle{amsplain}

\begin{thebibliography}{10} 

\bibitem{Andersen85} 
H.H. Andersen, 
\textit{Schubert varieties and Demazure's character formula}, 
Invent. Math. {\bf 79} (1985) 611--618. 
 
\bibitem {AntoineSpeiser64} 
J.-P. Antoine, D. Speiser,  
\textit{Characters of irreducible representations of the simple groups. 
I. General theory.}  
J. Math. Phys. {\bf 5} (1964) 1226--1234;
\hfill\break 
\textit{Characters of irreducible representations of the simple groups. II. 
Application to classical groups},  
J. Math. Phys. {\bf 5} (1964) 1560--1572.

\bibitem{Brion88} 
M. Brion, 
\textit{Points entiers dans les poly\`edres convexes}, 
Ann. scient. \'Ec. Norm. Sup., 4e s\'erie {\bf 21} (1988) 653--663.

\bibitem {Brion92} 
M. Brion, 
\textit{Poly\`edres et r\'eseaux}, 
Enseign. Math. (2) \textbf{38} (1992) 71--88. 

\bibitem{Demazure74} 
M. Demazure, 
\textit{D\'esingularisation des vari\'et\'es de Schubert g\'en\'eralis\'ees}, 
Annales scientifiques de l'\'E.N.S. $4^{\text e}$ s\'erie, tome 7 (1974) 53--88;
\hfill\break
\textit{Une nouvelle formule des caract\`eres}, 
Bull. Sci. Math. {\bf 98} (1974) 163--172.

\bibitem{DhillonKhare16} 
G. Dhillon, A. Khare, 
\textit{Characters of highest weight modules and integrability}, 
arXiv:1606.09640 (2016);
\hfill\break  
\textit{The Weyl-Kac weight formula}, 
S\'eminaire Lotharingien de Combinatoire \textbf{78B} (2017), 
Proceedings of the 29th Conference on Formal Power
Series and Algebraic Combinatorics (London), Article \#{77},  
arXiv:1802.06974 (2018).
 
\bibitem{Joseph85} 
A. Joseph, 
\textit{On the Demazure character formula}, 
Annales scientifiques de l'\'E.N.S. $4^{\text e}$ s\'erie, tome 18 (1985) 389--419.
 
\bibitem{Kass91} 
S. Kass, 
\textit{A recursive formula for characters of simple Lie algebras}, 
J. Alg. {\bf 137} (1991) 126--144. 

\bibitem{LecouveyLenart21}
C. Lecouvey, C. Lenart, 
\textit{Atomic decomposition of characters and crystals},
Adv. Math. {\bf 376} (2021) 107453. 

\bibitem{Macdonald72} 
I.G. Macdonald, 
\textit{Poincar\'e series of a Coxeter group}, 
Math. Ann. {\bf 199} (1972) 161--174.

\bibitem{Postnikov09} 
A. Postnikov, 
\textit{Permutohedra, associahedra, and beyond}, 
Int. Math. Res. Not. IMRN {\bf 6} (2009) 1026--1106.

\bibitem{Rasmussen18} 
J. Rasmussen,  
 \textit{Layer structure of irreducible Lie algebra modules},
 arXiv:1803.06592 (2018).

\bibitem{Schutzer12} 
W. Schutzer, 
\textit{A new character formula for Lie algebras and Lie groups}, 
J. Lie Theory 22.3 (2012) 817--838. 

\bibitem{Walton04} 
M.A. Walton, 
\textit{Polytope sums and Lie characters}, 
Symmetry in physics, CRM Proc. Lecture Notes {\bf 34} (2004) 203--214, 
Proceedings of a CRM Workshop held in Memory of Robert T. Sharp, 12-14 September 2002. 

\bibitem{Walton13} 
M.A. Walton, 
\textit{Polytope expansion of Lie characters and applications}, 
J. Math. Phys. {\bf 54} (2013) 121701.  

\bibitem{MWproc} 
M.A. Walton, 
\textit{Demazure formulas for weight polytopes}, in\,\ 
M.B. Paranjape et al (eds.), Quantum Theory and Symmetries; 
Proceedings of XIth International Symposium, Montr\'eal, 1-5 July 2019 (Springer, 2021), pages 287--296.

\end{thebibliography}

\vskip1cm\noindent{\underline{Data availability statement}}\hfb 
Data sharing is not applicable to this article as no new data were created or analyzed in this
study.

\end{document}